\documentclass[12pt,a4paper]{amsart}

\theoremstyle{plain}
\newtheorem{lem}{Lemma}

\theoremstyle{definition}
\newtheorem{obs}{Observation}

\newcommand{\itr}[2]{#1^{\circ #2}}

\newcommand{\NN}{\mathbb{N}}
\newcommand{\ie}{i.\,e. }
\newcommand{\Aut}{GA_n}
\newcommand{\LF}{LF_n}
\newcommand{\F}{\mathcal{F}_n}

\newcommand{\Diag}{Di_n}
\newcommand{\Aff}{Af_n}
\newcommand{\Elem}{El_n}
\newcommand{\Tame}{T_n}
\newcommand{\Jac}{J_n}

\title{Remarks on a normal subgroup of $GA_n$}
\author[J. Zygad{\l}o]{Jakub Zygad{\l}o}
\date{April 9, 2008}
\subjclass[2000]{14R10; 20E99}
\address{Jagiellonian University \\
Institute of Mathematics \\
Reymonta 4 \\
30-059 Krak\'ow \\
Poland.}
\email{jakub.zygadlo@im.uj.edu.pl}

\begin{document}
\begin{abstract}
We show that the subgroup generated by locally finite polynomial automorphisms of $k^n$ is normal in $GA_n$. Also, some properties of normal subgroups of $GA_n$ containing all diagonal automorphisms are given.
\end{abstract}
\maketitle 

Let $k$ be a field of characteristic zero. A {\it polynomial endomorphism} $G$ of $k^n$ will be identified with an $n$-tuple $G=(G_1,\ldots,G_n)$ of polynomials in $n$ variables over $k$. Every polynomial endomorphism of $k^n$ defines a $k$-endomorphism $G^*$ of the algebra $k[X_1,\ldots,X_n]$ in the following way: $G^*(f):=f\circ G=f(G_1,\ldots,G_n)$ (equivalently, set $G^*(X_i):=G_i$, for $i=1,\ldots,n$ and extend $G^*$ uniquely to a $k$-endomorphism). A polynomial endomorphism $G$ of $k^n$ is an {\it automorphism}, if its associated map $G^*$ is a $k$-automorphism of the algebra $k[X_1,\ldots,X_n]$. The set of all polynomial automorphisms of $k^n$ forms a group under composition - it is called the {\it polynomial automorphism group of $k^n$} and denoted $\Aut$.\\
Given an endomorphism $G=(G_1,\ldots,G_n)$ of $k^n$, we will mean by $\deg G_i$ the total degree of $G_i$ and $\deg G:=\max_{i=1,\ldots n}\deg G_i$. Notation $\itr{G}{m}$ will denote the $m$-th iterate of $G$, \ie $\itr{G}{m}=G\circ\itr{G}{(m-1)}$ for $m\geq 1$ and $\itr{G}{0}:=I$ - the identity.

As in \cite{FM}, an automorphism $G\in\Aut$ will be called {\it locally finite} if $G^*$ is locally finite as a $k$-linear mapping - that is, for every $f\in k[X_1,\ldots,X_n]$ the $k$-vector subspace of $k[X_1,\ldots,X_n]$ spanned by 
$$\{\itr{(G^*)}{m}(f):m\in\NN\}=\{f\circ\itr{G}{m}:m\in\NN\}$$
is finite dimensional. The set of all locally finite polynomial automorphisms of $k^n$ will be denoted $\LF$. It is not hard to see that the following conditions on $G$ are equivalent (see \cite{FM}, Th. 1.1):
\begin{enumerate}
\item $G$ is locally finite,
\item $\sup_{m\in\NN}\deg(\itr{G}{m})<+\infty$,
\item $\exists p(T)\in k[T]\setminus\{0\}: p(G)=0$, that is $a_d\itr{G}{d}+\ldots+a_1G+a_0I=0$ for some $a_i\in k$, $i=0,\ldots,d$, not all $a_i$ equal 0.
\end{enumerate}
Any $p$ satisfying condition $(3)$ is called the {\it vanishing polynomial} for G. One easily verifies that the set of all vanishing polynomials (for a given $G$) together with the zero polynomial forms an ideal in $k[T]$ - its monic generator will be called the {\it minimal polynomial} for $G$ and denoted by $\mu_G$.

Fix $n\geq 1$ and consider the subgroup $\F\subset \Aut$ generated by all locally finite polynomial automorphisms of $k^n$:
$$\F=\{G_1\circ\ldots\circ G_s: G_i\in\LF\text{ for }i=1,\ldots,s\text{ and }s\in\NN\}$$

\begin{obs}
$\F$ is a normal subgroup of $\Aut$.
\begin{proof}
Obviously $\F$ contains $I$ and is closed under composition. Since for $F=G_1\circ\ldots\circ G_s$ and $\phi\in\Aut$ we have $F^{-1}=G_s^{-1}\circ\ldots\circ G_1^{-1}$ and $\phi\circ F\circ\phi^{-1}=(\phi\circ G_1\circ\phi^{-1})\circ\ldots\circ(\phi\circ G_s\circ\phi^{-1})$, it remains to show that if $G$ is locally finite, then so are $G^{-1}$ and $\phi\circ G\circ\phi^{-1}$. But if $p(T)=a_dT^d+\ldots+a_1T+a_0$ is a vanishing polynomial for $G$, then $T^dp(1/T)=a_0T^d+a_1T^{d-1}+\ldots+a_d$ vanishes on $G^{-1}$. Also note that $\sup_{m\in\NN}\{\deg\itr{(\phi\circ G\circ\phi^{-1})}{m}\}= \sup_{m\in\NN}\{\deg(\phi\circ\itr{G}{m}\circ\phi^{-1})\}\leq \deg\phi\cdot\sup_{m\in\NN}\{\deg(\itr{G}{m})\}\cdot\deg(\phi^{-1})<+\infty$ if $G\in\LF$. This concludes the proof.
\end{proof}
\end{obs}

An automorphism $G=(G_1,\ldots,G_n)\in\Aut$ is called
\begin{itemize}
\item {\it diagonal}, if $G_i=c_iX_i$, $c_i\in k\setminus\{0\}$ for $i=1,\ldots,n$,
\item {\it affine}, if $G_i=\sum_{j=1}^na_{ij}X_j+b_i$, $a_{ij}\in k$, $b_i\in k$ for $i,j=1,\ldots,n$
\item {\it elementary}, if $G_i=X_i+g$ for some $g\in k[X_1,...,X_{i-1},X_{i+1},...,X_n]$
and $G_j=X_j$ for $j\neq i$,
\item {\it tame}, if $G$ is a (finite) composition of affine and elementary automorphisms,
\item {\it jacobian}, if $\det J(G)=1$, where $J(G)=\big(\frac{\partial G_i}{\partial X_j}\big)_{i,j=1,\ldots,n}$ is the jacobian matrix of $G$
\end{itemize}
The set of diagonal (resp. affine, elementary, tame, jacobian) automorphisms of $k^n$ will be denoted by $\Diag$ (resp. $\Aff$, $\Elem$, $\Tame$, $\Jac$). We will need the following standard lemma:

\begin{lem}\label{tames}
Every tame automorphism $F\in\Tame$ can be written in the form
$$F=E_1\circ\ldots\circ E_s\circ D$$
where $E_i\in\Elem$ for $i=1,\ldots,s$ and $D\in\Diag$.
\begin{proof}
First, note that if $E=(X_1,\ldots,X_{i-1},X_i+g,X_{i+1},\ldots,X_n)\in\Elem$ and $D\in\Diag$ then $D\circ E=\tilde{E}\circ D$, where $\tilde{E}=(X_1,\ldots,X_i+\tilde{g},\ldots,X_n)$ and $\tilde{g}=g\circ D^{-1}=g(D^{-1}_1,\ldots,D^{-1}_{i-1},D^{-1}_{i+1},\ldots,D^{-1}_n)$. Therefore it suffices to prove that every $F\in\Tame$ is a composition of elementary and diagonal maps, and consequently that the claim holds for $F\in\Aff$. Since translations (maps of the form $B_i=X_i+b_i$, $b_i\in k$, $i=1,\ldots,n$) are clearly compositions of elementary morphisms,
we may assume that $F(0)=0$. Consider the matrix $A$ of $F$ in some basis of $k^n$; of course $\det A\neq 0$. Applying Gaussian elimination to $A$, we arrive at a triangular matrix $A'$ associated to a mapping $F'$ that clearly is a composition of elementary automorphisms. To conclude the proof, we need to show that the following automorphisms (operations in the procedure of Gaussian elimination) can be expressed as compositions of elementary and diagonal maps:
\begin{itemize}
\item[-] $A_{i,j,c}=(X_1,\ldots,X_{i-1},X_i+cX_j,X_{i+1},\ldots,X_n)$ - addition of a scalar multiple of a row to another row,
\item[-] $T_{i,j}=(X_1,\ldots,X_{i-1},X_j,X_{i+1},\ldots,X_{j-1},X_i,X_{j+1},\ldots,X_n)$ - transposition of rows.
\end{itemize}
But clearly $A_{i,j,c}$ is elementary and $T_{i,j}=A_{i,j,1}\circ A_{j,i,-1} \circ A_{i,j,1}\circ\tilde{D}$, where $\tilde{D}:=(X_1,\ldots,X_{i-1},-X_i,X_{i+1},\ldots,X_n)\in\Diag$.
\end{proof}
\end{lem}

Here is our main observation:

\begin{obs}\label{nor}
Let $N$ be a normal subgroup of $\Aut$ containing all diagonal automorphisms. Then every tame automorphism belongs to $N$.
\begin{proof}
Due to Lemma \ref{tames}, we only need to prove that every elementary automorphism is in $N$. Let $F\in\Elem$, that is
$$F=(X_1,\ldots,X_{i-1},X_i+g,X_{i+1},\ldots,X_n)$$
for some $i\in\{1,\ldots,n\}$ and $g\in k[X_1,\ldots,X_{i-1},X_{i+1},\ldots,X_n]$. Take $D:=(X_1,\ldots,X_{i-1},2X_i,X_{i+1},\ldots,X_n)$ and note that
$$F^{-1}\circ D\circ F=(X_1,\ldots,X_{i-1},2X_i+g,X_{i+1},\ldots,X_n)\in N$$ since $N$ is normal. Then $(F^{-1}\circ D\circ F)\circ D^{-1}=F$ is an element of $N$.
\end{proof}
\end{obs}

The following observation (applied to the case $n=2$) can be of some revelance to the results of \cite{Da}.

\begin{obs}\label{unip}
Let $N$ be a normal subgroup of $\Jac$, such that every diagonal automorphism with determinant equal to 1 belongs to $N$. Then $\Tame\cap\Jac\subset N$.
\begin{proof}
Express $F\in\Tame\cap\Jac$ as in Lemma \ref{tames} - since all $E_i\in\Jac$, we conclude that $D\in\Jac$ and as in the previous observation, it suffices to show that $\Elem\subset N$. Let $F=(X_1,\ldots,X_{i-1},X_i+g,X_{i+1},\ldots,X_n)\in\Elem$, with $g\in k[X_1,\ldots,\hat{X_i},\ldots,X_n]$. Fix any $j\neq i$ and write $g=\sum_{r=0}^dg_r(X_*)X_j^r$, where $X_*$ denotes the set of variables with $X_i$ and $X_j$ omitted. Next, choose $a\in k\setminus\{0\}$ that is not a root of unity (such an $a$ exists, since $k$ is of characteristic 0) and define mappings: $D\in\Diag\cap\Jac$ as $D_i=aX_i$, $D_j=a^{-1}X_j$ and $D_l=X_l$ for $l\not\in\{i,j\}$ and $E\in\Elem$ as $E_i=X_i+h$, where $h=\sum_{r=0}^d(a^{1+r}-1)^{-1}g_r(X_*)X_j^r$ and $E_l=X_l$ for $l\neq i$. Then $E^{-1}\circ D\circ E\in N$ and
$$(E^{-1}\circ D\circ E)_l=\begin{cases}
aX_i+ah-h\circ D, &l=i\\
a^{-1}X_j, &l=j\\
X_l, &l\not\in\{i,j\}
\end{cases}$$
It follows that $(E^{-1}\circ D\circ E\circ D^{-1})_i=X_i+ah\circ D^{-1}-h$ (identity on other coordinates). But $ah\circ D^{-1}-h=\sum_{r=0}^d(a^{1+r}-1)h_r(X_*)X_j^r=g$, hence $F=(E^{-1}\circ D\circ E)\circ D^{-1}\in N$.
\end{proof}
\end{obs}

Surprisingly, similar fact occurs for the famous Nagata automorphism (see \cite{Na}).

\begin{obs}\label{Nag}
Let $n=3$ and $N$ be a normal subgroup of $GA_3$ containing all diagonal automorphisms. Then the Nagata automorphism
$$F=(X-2(Y^2+XZ)Y-(Y^2+XZ)^2Z,Y+(Y^2+XZ)Z,Z)$$
belongs to $N$.
\begin{proof}
Let $\sigma:=Y^2+XZ$. Take $L:=(\frac{1}{4}X,\frac{1}{2}Y,Z)\in N$ and note that $\sigma\circ L=\frac{1}{4}\sigma$. One easily checks that $F^{-1}=(X+2\sigma Y-\sigma^2 Z,Y-\sigma Z,Z)$ and
$$F^{-1}\circ L\circ F=\big(\frac{1}{4}X-\frac{1}{4}\sigma Y-\frac{1}{16}\sigma^2 Z,\frac{1}{2}Y+\frac{1}{4}\sigma Z,Z\big)$$
Consequently, $(F^{-1}\circ L\circ F)\circ L^{-1}=F$ is an element of $N$. 
\end{proof}
\end{obs}

Since $\F$ is normal for every $n\in\NN$ and all diagonal automorphisms are clearly locally finite, we have the following alternative:
\begin{enumerate}
\item either $\F$ is a nontrivial normal subgroup of $\Aut$,
\item or $\F=\Aut$, and thus locally finite automorphisms form the set of generators for $\Aut$.
\end{enumerate}
Note that for $n=1,2$ the second possibility occurs (case $n=1$ being trivial and $n=2$ a consequence of Jung-van der Kulk theorem and Observation \ref{nor}). For $n=3$, theorem of Shestakov-Umirbaev and Observation \ref{Nag} may suggest that $\mathcal{F}_3=GA_3$ but the problem is certainly worth further research.

\end{document}